\documentclass{amsart}

\usepackage{amsfonts,amsthm,amssymb}
\usepackage{epsfig}
\usepackage{pstricks,slashbox,multirow}
\theoremstyle{plain}
\newtheorem{theorem}{Theorem}[section]

\newtheorem{lemma}[theorem]{Lemma}
\newtheorem{prop}[theorem]{Proposition}

\newtheorem{obs}[theorem]{Observation}

\newtheorem{definition}[theorem]{Definition}

\newtheorem{exm}[theorem]{Example}

\newtheorem{rem}[theorem]{Remark}

\begin{document}

\date{}

\title[Generating Primes - A New Approach]{Generating Primes Numbers - A Fast New Approach}

\author{V.Vilfred} 

\address{Formerly Professor, Department of Mathematics, Central University of Kerala, Tejaswini Hills, Periye - 671 316, Kasaragod \\ Kerala, India.}

\email{vilfredkamal@gmail.com}


\begin{abstract}
Bertrand's Postulate ensures existence of prime $p$ between $n$ and $2n$, $n$ an integer $\geq 2$ and the sieve of Eratosthenes, a very simple ancient algorithm, generates all prime numbers up to any given limit. Combining the above two, in this paper, we provide a simple fast moving algorithm to generate prime numbers up to any given limit.  We also discuss Riemann zeta function related to generating of prime numbers. 
\end{abstract}

\maketitle


------------------

\section{Introduction} 

In generating prime numbers, the sieve of Eratosthenes \cite{dg16,mmp16} and Bertrand postulate \cite{b45} are important mile stones. The sieve of Eratosthenes is a very simple ancient algorithm that generates all primes up to any given limit. In 1845, Joseph Bertrand postulated that for any integer $n > 3$, there exists prime $p$ such that $n < p < 2n - 2$. And its slightly weaker form is that there exists prime $p$ such that $n < p < 2n$, $n$ an integer $\geq 2$. In 1850, Pafnuty Chebychev \cite{t1852} first proved this postulate analytically. In 1932, Paul Erdos \cite{dg16} gave an elementary proof using facts about the middle binomial coefficient. In 2002, Manindra Agarwal, Neeraj Kayal and Nitin Saxena \cite{aks04} presented an unconditional deterministic polynomial-time algorithm that determines whether an input number is prime or composite.

The problem of computing $\pi (x)$, the number of primes less than or equal to $x$ is one of the oldest problem in Mathematics, $x\in \mathbb{N}$. For a  very long time, the sieve of Eratosthenes has been the practical way to compute $\pi(x)$ despite its time complexity. Legendre \cite{l59} observed  a  combinatorial formula, known as {\em Legendre sum}, for the number of primes $p$ for which $x^{\frac{1}{2}} < p \leq x$. Since then, a large number  of  writers have  suggested variants and improvements of the formula. During 1870 to 1885, astronomer Meissel \cite{m85} developed practical combinatorial method to compute $\pi(x)$ and in  1959, Lehmer \cite{l59} extended  and simplified  Meissel's  method. In 1985, the Meissel-Lehmer method was used to compute several values of $\pi(x)$ up to $x = 4.10^{16}$ \cite{lmo} and in 1996, Deleglise and Rivat \cite{dr} developed a modified form of the Meissel-Lehmer method saving much computation.       

In this paper, using both the sieve of Eratosthenes and Bertrand postulate, we obtain prime numbers in a faster way by considering $[p_i+1, p_i^2+4p_i+3]$ as consecutive intervals where $p_1, p_2, \ldots$ denotes the primes 2, 3, \ldots numbered in increasing order, $i \geq 2$. It is a tight bound in the sense that if we increase the width of the interval still further, then the algorithm will not work, in general. We present two more methods to show the fastness of our method. In all the three methods, we obtain primes in each successive intervals of the form $[p_i+1, 2p_j]$ and existence of prime(s) in each such intervals is ensured by Bertrand postulate, $j \geq i \geq 2$. In the first method, we consider successive intervals $[p_i+1, 2p_i]$, $i \geq 2$. In the second method, successive intervals $[p_i+1, p_i^2]$ are considered, $i \geq 2$. In the third method, we consider $[p_i+1, p_i^2+4p_i+3]$ as consecutive intervals. The author feels that this development is going to revolutionise development in Mathematics, especially in Cryptography, Number theory, Signal processing and Computational Mathematics.    

Even though the third method generates primes in a faster way and is the best, we present the other two to highlight their differences. After generating, it is easy to check, from the listing, whether any given number is prime or not. We also present Riemann zeta function \cite{jps} related to prime generation.   

  Algorithms used in the three methods to generate prime numbers up to any given limit are 
 \begin{enumerate} 
 	\item [\rm (i)]  $p_1 = 2$, $IpP_1 = \{2\}$ = $pP_1$, $MIpP_1$ = 2 = $MpP_1$, $\#pP_1$ = 1 and other consecutive intervals $IpP_{j+1}$ = $[MIpP_j+1, 2MpP_j]$ and $pP_{j+1}$ = $\varphi(IpP_{j+1})$ = set of primes in $IpP_{j+1}$ for $j = 1, 2, \ldots$; 
   \item [\rm (ii)]  $p_1$ = 2, $IpP_1$ = $\{2\}$ = $pP_1$, $MIpP_1$ = 2 = $MpP_1$, $\#pP_1$ = 1, $p_2$ = 3, $IpP_2$ = $\{3\}$ = $pP_2$, $MIpP_2$ = 3 = $MpP_2$, $\#pP_2$ = 1, $IpP_{j+1}$ = $[MIpP_{j}+1, MpP_j^2]$ and $pP_{j+1}$ = $\varphi([MIpP_{j}+2, MpP_j^2-2])$ for $j = 2, 3, \ldots$ and 
    \item [\rm (iii)]  $p_1$ = 2, $p_2$ = 3, $IpP_1$ = $[2, 3]$, $pP_1$ = $\{2, 3\}$, $MIpP_1$ = 3 = $MpP_1$, $\#pP_1$ = 2, $MIpP_{j}$ = max $IpP_{j}$, $MpP_{j}$ = max $pP_{j}$, $\#pP_j$ = $|pP_j|$, $IpP_{j+1}$ = $[MIpP_{j}+1, MpP_j^2 +4MpP_j+3]$ and $pP_{j+1}$ = $\varphi(IpP_{j+1})$ = $\varphi([MIpP_{j}+1, MpP_j^2+4MpP_j+2])$ for $j$ = $1, 2, \ldots$.
    
    In all the three methods, $MIpP_{j}$ = max $IpP_{j}$ = maximum value of $IpP_{j}$, $MpP_{j}$ = max $pP_{j}$ = maximum value in $pP_{j}$, $\#pP_j$ = $|pP_j|$ and $pP_{j+1}$ = $\varphi(IpP_{j+1})$ = set of primes in $IpP_{j+1}$ for $j$ = $1, 2, \ldots$.
\end{enumerate} 
 
 Effort to find a faster method to generate prime numbers compared with  
 the sieve of Eratosthenes is the motivation for this work. For all basic ideas in number theory, we follow \cite{tk02}. Throughout the paper, for $m,n\in \mathbb{N}$ and $m \leq n$, $[m, n]$ $= \{k\in \mathbb{N}: m \leq k \leq n\}$ and $p_1, p_2, \ldots$ denotes the primes 2, 3, \ldots numbered in increasing order.

\section{Preliminaries}   

To simplify our work, the following notations are used in this paper. 

$\mathbb{N}$ = $\{1,2,\ldots\}$; $\mathbb{N}_0$ = $\mathbb{N} \cup \{0\}$ = $\{0,1,2,\ldots\}$;

$\mathbb{P}$ = the set of all prime numbers = $\varphi(\mathbb{N})$; 

$\mathbb{C}$ = the set of all composite numbers so that 

$\mathbb{P} \cap \mathbb{C}$ = $\emptyset$ and  $\mathbb{P} \cup \mathbb{C}$ = $\mathbb{N} \setminus \{1\}$; 

$\mathbb{P}(S) = \varphi(S) =$ the set of all primes in $S$ and 

$\mathbb{C}(S)$ = the set of all composite numbers in $S$, $S \subseteq \mathbb{N}$;

$\pi(S) =$ the number of primes in $S$ when $S$ is finite and $S \subset \mathbb{N}$; 

$\varphi(n)$ = $\varphi([1, n])$ = the set of all primes $\leq n$ and

$\pi(n)$ = $\pi([1, n])$ = the number of primes $\leq n$, $n\in \mathbb{N}$;

$\left\lfloor \frac{m}{n}\right\rfloor$ = integer part of $\frac{m}{n}$ when $m,n\in\mathbb{N}$.
\vspace{.2cm}

\begin{definition} Let $n_1 \leq n_2$, $a,n_1,n_2\in \mathbb{N}$ and $[a]$ $= \{a, 2a, \ldots \} = a\mathbb{N}$, the set of multiples of $a$ in $\mathbb{N}$. Then, we denote the set of all multiples of $a$ each lies between $n_1$ and $n_2$ by $([a]: n_1, n_2)$. Thus, $([a]: n_1, n_2)$ $= \{ka/~ n_1 \leq ka \leq n_2, ~ k\in \mathbb{N}\}$, $a,n_1,n_2\in \mathbb{N}$, $n_1 \leq n_2$.
\end{definition}

Consider the following simple lemma which is an important result used in this paper to generate larger primes.

\begin{lemma} \label{b1} {\rm Let $1 \leq i < j$, $i,j\in \mathbb{N}$, $p_i, p_j$ be primes, $Q_{i,j}$ = $\left\lfloor \frac{p_j}{p_i} \right\rfloor$, $Q_{i,j}^{'}$ = $\left\lfloor \frac{2p_j}{p_i}\right\rfloor$, $Q_{i,j}^{''}$ = $\left\lfloor \frac{p_j^2}{p_i} \right\rfloor$ and $Q_{i,j}^{'''}$ = $\left\lfloor \frac{p_j^2+4p_j+3}{p_i} \right\rfloor$. Then, 
\begin{enumerate} 
\item [\rm (a)] \label{l11} $p_i(Q_{i,j}+1)$ is the smallest integer multiple of $p_i$ that is greater than $p_j$.
\item [\rm (b)]  \label{l12} $p_iQ_{i,j}^{'}$ is the biggest integer multiple of $p_i$ that is less than or equal to $2p_j$.
\item [\rm (c)]  \label{l13} $p_iQ_{i,j}^{''}$ is the biggest integer multiple of $p_i$ that is less than or equal to $p_j^2$.
\item [\rm (d)] \label{l14} $p_iQ_{i,j}^{'''}$ is the biggest integer multiple of $p_i$ that is less than or equal to $p_j^2+4p_j+3$.
\item [\rm (e)] \label{l15} $([p_i]:p_j+1, 2p_j)$ $= \{p_i(Q_{i,j}+1), p_i(Q_{i,j}+2)$, $\ldots, p_iQ_{i,j}^{'}\}$.
\item [\rm (f)] \label{l16} $([p_i]:p_j+1, p_j^2)$ $= \{p_i(Q_{i,j}+1), p_i(Q_{i,j}+2), \ldots, p_iQ_{i,j}^{''}\}$. 
\item [\rm (g)] \label{l17} $([p_i]:p_j+1, p_j^2+4p_j+3)$ $= \{p_i(Q_{i,j}+1), p_i(Q_{i,j}+2), \ldots, p_iQ_{i,j}^{'''}\}$.
\item [\rm (h)] \label{l18} In $[p_j+1, p_j^2+4p_j+3]$, any composite number has $p_1$, $p_2$, $\ldots$, $p_{j-1}$ or $p_j$ as a divisor. And $p_j^2+4p_j+3 < p_{j+1}^2$ for $j \geq 2$.
\item [\rm (i)] \label{l19} For $j \geq 2$, any composite number whose prime divisors, each $> p_j$ will be $\geq p_{j+1}^2 \geq {(p_j+2)}^2$ $>$ $p_j^2+4p_j+3$.
\end{enumerate} }
\end{lemma}

\begin{proof} Given, $1 \leq i < j$, $Q_{i,j} = \left\lfloor \frac{p_j}{p_i} \right\rfloor$, $Q_{i,j}^{'} = \left\lfloor \frac{2p_j}{p_i}\right\rfloor$, $Q_{i,j}^{''} = \left\lfloor \frac{p_j^2}{p_i} \right\rfloor$ and $Q_{i,j}^{'''} = \left\lfloor \frac{p_j^2+4p_j+3}{p_i} \right\rfloor$. This implies, $p_i < p_j$ and $Q_{i,j}$, $Q_{i,j}^{'}$, $Q_{i,j}^{''}$, $Q_{i,j}^{'''}$ are the quotients when $p_j$, $2p_j$, $p_j^2$, $p_j^2+4p_j+3$ are divided by $p_i$, respectively. 

Let $p_j = p_iQ_{i,j} + R_{i,j}$ where $Q_{i,j}$ and $R_{i,j}$ are quotient and remainder when $p_j$ is divided by $p_i$, $1 \leq R_{i,j} \leq p_i -1$ since $p_j$ is a prime number greater than $p_i$. Similarly, let $2p_j = p_iQ_{i,j}^{'} + R_{i,j}^{'}$, $p_j^2 = p_iQ_{i,j}^{''} + R_{i,j}^{''}$ and $p_j^2+4p_j+3 = p_iQ_{i,j}^{'''} + R_{i,j}^{'''}$ where $Q_{i,j}^{'}$ and $R_{i,j}^{'}$ are quotient and remainder when $2p_j$ is divided by $p_i$, $Q_{i,j}^{''}$ and $R_{i,j}^{''}$ are quotient and remainder when $p_j^2$ is divided by $p_i$ and $Q_{i,j}^{'''}$ and $R_{i,j}^{'''}$ are quotient and remainder when $p_j^2+4p_j+3$ is divided by $p_i$,  $0 \leq R_{i,j}^{'}, R_{i, j}^{'''} \leq p_i -1$ and $1 \leq  R_{i, j}^{''} \leq p_i -1$. From the above, we get results (a), (b), (c) and (d).

Result (e) follows from (a) and (b). Result (f) follows from (a), (c) and $p_j$ is an odd prime $> 2$. Result (g) follows from (a) and (d).

For $j \geq 2$, $p_{j+1} \geq p_j+2 > p_j$. This implies, the smallest composite number whose divisors are other than $p_1$, $p_2$, $\ldots$, $p_{j-1}$ and $p_j$ is $p_{j+1}^2 \geq {(p_j+2)}^2$ = $p_j^2+4p_j+4 > p_j^2+4p_j+3$ for $j \geq 2$. Hence, we get results (h) and (i). 
\end{proof}

\begin{rem} \quad
	In the above lemma, results $(h)$ and $(i)$ are not true for $j$ = 1 since $p_2^2$ = $3^2$ = 9 $<$ $p_1^2+4p_1+3$ = $2^2+8+3$ = 15 $<$ 16 = $p_1^2+4p_1+4$ = ${(p_1+2)}^2$. 
\end{rem}

\section{Main Result}

In this section, we present all the three methods to generate prime numbers even though the third is the best. All the three methods use the sieve of Eratosthenes \cite{dg16,mmp16}, Bertrand's postulate on prime numbers \cite{b45} and Lemma \ref{b1}. After generating primes, it is easy to check, from the listing whether any given number is prime or not.

In all the three methods, using Lemma \ref{b1} we calculate all composite numbers, then all prime numbers in the interval $[p_j+1, 2p_j]$ using Theorems \ref{c2} and \ref{c3} in the first  method, in $[p_j+1, p_j^2]$ using Theorems \ref{c4} and \ref{c5} in the second method and finally in $[p_j+1, p_j^2+4p_j+3]$ using Theorems \ref{c6} and \ref{c7} in the third method, $j \geq 2$.     

\begin{lemma} \label{c1} {\rm For $1\leq i < j$ and $i,j\in \mathbb{N}$, 
\begin{enumerate} 
\item [\rm (a)] \label{l21} the set of all composite numbers in $[p_j+1, 2p_j]$ is given by 

$\mathbb{C}([p_j+1, 2p_j])$ $=$ $\bigcup^{j}_{i=1}$ $p_i\{Q_{i,j}+1,Q_{i,j}+2,\ldots,Q_{i,j}^{'}\}$ 

\hspace{2.3cm} = $\bigcup^{j-1}_{i=1}$ $p_i\{Q_{i,j}+1,Q_{i,j}+2,\ldots,Q_{i,j}^{'}\}$ $\bigcup$ $\{2p_j\}$;
\item [\rm (b)]  \label{l22} the number of composite numbers in $[p_j+1, 2p_j]$ is given by

 $\#\mathbb{C}([p_j+1, 2p_j])$ $= \sum_{i=1}^j({Q_{i,j}}^{'}-Q_{i,j})$; and
\item [\rm (c)] \label{l23} the number of primes in $[p_j+1, 2p_j]$ is given by

$\pi([p_j+1, 2p_j])$ $= 2p_2 - p_j - \#\mathbb{C}([p_j+1, 2p_j])$.
\end{enumerate} }
\end{lemma}
\begin{proof} For $1\leq i < j$, we have $Q_{i,j} = \left\lfloor \frac{p_j}{p_i} \right\rfloor$ and $Q_{i,j}^{'} = \left\lfloor \frac{2p_j}{p_i} \right\rfloor$. 
\begin{enumerate} 
\item [\rm (a)] Then, the set of all composite numbers in $[p_j+1, 2p_j]$ is
 
$\mathbb{C}([p_j+1, 2p_j])$ = $\bigcup^{j}_{i=1} \{$multiples of $p_i$ in $[p_j+1, 2p_j]\}$.

\hspace{2.7cm} = $\bigcup^{j}_{i=1}$ $([p_i]: p_j+1, 2p_j)$.

\hspace{.2cm} = $\{2p_j\}$  $\bigcup$ $(\bigcup^{j-1}_{i=1}$ $\{p_i(Q_{i,j}+1),p_i(Q_{i,j}+2),\ldots,p_i Q_{i,j}^{'} \})$ 

\hspace{.5cm} using $(e)$ of Lemma \ref{b1} and  also $2p_j$ is the only 

\hspace{.5cm} composite number with $p_j$ as a divisor in $[p_j+1, 2p_j]$. 
\item [\rm (b)] Composite numbers that are divisible by $p_i$ in $[p_j+1, 2p_j]$ are

$\{p_i(Q_{i,j}+1),p_i(Q_{i,j}+2),\ldots,p_i Q_{i,j}^{'} \}$, $1 \leq i \leq j$. 

 This implies, the number of composite numbers that are divisible by $p_i$ in $[p_j+1, 2p_j]$ is ${Q_{i,j}}^{'}-Q_{i,j}$. This implies, the number of composite numbers in $[p_j+1, 2p_j]$ is 
 
 $\#\mathbb{C}([p_j+1, 2p_j])$ $= \sum_{i=1}^j({Q_{i,j}}^{'}-Q_{i,j})$. 

\item [\rm (c)] Number of primes in $[p_j+1, 2p_j]$ 

= $|[p_j+1, 2p_j]|$ - number of composite numbers in $[p_j+1, 2p_j]$

= $\pi([p_j+1, 2p_j])$ $= 2p_2 - p_j - \#\mathbb{C}([p_j+1, 2p_j])$ using $(b)$.
\end{enumerate} 
Hence we get the result.
\end{proof}

{\subsection ~  {\bf Calculation of Prime numbers in $[p_j+1, 2p_j]$, $j\geq 2$}}

Here, we obtain all prime numbers contained in the interval $[p_j+1, 2p_j]$ by removing its composite numbers using Lemma \ref{c1}(a) for $j\geq 2$.

\begin{theorem} \label{c2} {\rm Let $1 \leq i <j$, $i,j\in \mathbb{N}$, $Q_{i,j} = \left\lfloor \frac{p_j}{p_i} \right\rfloor$, $Q_{i,j}^{'} = \left\lfloor \frac{2p_j}{p_i} \right\rfloor$ and $S_j$ = $\varphi([p_j+1, 2p_j])$ = the set of all primes in $[p_j+1, 2p_j]$. Then for $j \geq 2$,
\begin{enumerate}
	\item [\rm (a)] $S_j$ = $[p_j+1, 2p_j]  \setminus \mathbb{C}([p_j+1, 2p_j])$ = $[p_j+2, 2p_j-1] \setminus \mathbb{C}([p_j+2, 2p_j-1])$;
	\item [\rm (b)] $S_j$ is non-empty;
	\item [\rm (c)]  $S_j$ = $[p_j+1, 2p_j]$ $\setminus$ $(\bigcup^{j}_{i=1}$ $p_i\{Q_{i,j}+1,Q_{i,j}+2,\ldots,Q_{i,j}^{'}\})$; and
	\item [\rm (d)] $S_j$  = $[p_j+2, 2p_j-1]$ $\setminus$ $(\bigcup^{j-1}_{i=1}$ $p_i\{Q_{i,j}+1,Q_{i,j}+2,\ldots,Q_{i,j}^{'}\})$.	 
\end{enumerate}  }
\end{theorem}
\begin{proof} For $j \geq 2$, the set of all prime numbers contained in $[p_j+1, 2p_j]$ is obtained by removing all composite numbers contained in $[p_j+1, 2p_j]$. 

For $j \geq 1$, by  Bertrand postulate, $[p_j+1, 2p_j]$ contains at least one prime and for $j \geq 2$, $p_j+1$ is composite. Hence, results $(a)$ and $(b)$ are true.  

For $j \geq 2$, using $(1)$, $[p_j+2, 2p_j]$ contains at least one prime. And any prime number contained in $[p_j+2, 2p_j]$ is greater than $p_j$ and less than $2p_j$. But for $j \geq 2$, any composite number, in $[p_j+2, 2p_j]$, contains $p_1$, $p_2$, $\ldots$ or $p_j$ as a factor. Hence, after removing all multiples of $p_1$, $p_2$, $\ldots$, $p_j$ from $[p_j+2, 2p_j]$, the resultant set contains only prime(s). This implies, for $j \geq 2$,
\\ 
$S_j$ = $[p_j+1, 2p_j]$ $\setminus$ $\mathbb{C}([p_j+1, 2p_j])$ 

$= [p_j+2, 2p_j-1]$ $\setminus$ $\mathbb{C}([p_j+2, 2p_j-1])$ = $\varphi([p_j+2, 2p_j-1])$ 

= $[p_j+2, 2p_j-1]$ $\setminus$ ($\bigcup^{j-1}_{i=1}$ $p_i\{Q_{i,j}+1,Q_{i,j}+2,\ldots,Q_{i,j}^{'}\})$ using Lemma \ref{c1} where $Q_{i,j} = \left\lfloor \frac{p_j}{p_i} \right\rfloor$ and $Q_{i,j}^{'} = \left\lfloor \frac{2p_j}{p_i} \right\rfloor$, $1 \leq i < j$. And thereby $(c)$ and $(d)$ are true.
\end{proof}
In the above result, we could see that $\varphi([p_j+2, 2p_j-1])$, the set of all prime numbers in the interval $[p_j+1, 2p_j]$, is a non-empty set for every $j \geq 2$, $j\in \mathbb{N}$ where $p_j$ is the $j^{th}$ prime with $p_1 = 2$, $p_2 = 3$, $p_3 = 5$, $\ldots$. In Example \ref{e3}, using Theorem \ref{c2}, we calculate the number of prime numbers in $IpP_j$ = $\pi(IpP_j)$ = $\pi(pP_j) = |pP_j|$ with the notation of $MIpP_{j}$ = maximum value in $IpP_{j}$ and $MpP_{j}$ = maximum value in $pP_{j}$, $j$ = $1, 2, \ldots$. 

\begin{exm} \label{e3} $p_1 = 2$, $IpP_1 = \{ 2\}$ $= pP_1$, $MIpP_1 = 2 = MpP_1$. 

$\Rightarrow$ $IpP_2 = [MIpP_1+1, 2MpP_1] = [3, 4]$, $MIpP_2 = 4$,

\hspace{.5cm} $pP_2 = \varphi(IpP_2) = \varphi([3, 4]) = \{3\}$, $p_2 = 3$, $MpP_2 = 3$.

$\Rightarrow$ $IpP_3 = [MIpP_2+1, 2MpP_2] = [5, 6]$, $MIpP_3 = 6$,

\hspace{.5cm} $pP_3 = \varphi(IpP_3)  = \varphi([5, 6])$ = $\{ 5\}$, $p_3 = 5$, $MpP_3 = 5$.

$\Rightarrow$  $IpP_4 = [MIpP_3+1, 2MpP_3]$ = $[7, 10]$, $MIpP_4 = 10$, 

\hspace{.5cm} $pP_4 = \varphi(IpP_4) = \{ 7\}$, $p_4 = 7$, $MpP_4 = 7$.

$\Rightarrow$  $IpP_5 = [MIpP_4+1, 2MpP_4] = [11, 14]$, $MIpP_5 = 14$,  

\hspace{.5cm} $pP_5 =$ $\varphi(IpP_5) = \varphi([11, 14])$ = $\{ 11, 13\}$, 

\hspace{.5cm} $MpP_5 = 13$, $p_5 = 11$, $p_6 = 13$.

$\Rightarrow$  $IpP_6 = [MIpP_5+1, 2MpP_5] = [15, 26]$, $MIpP_6 = 26$,  

\hspace{.6cm} $pP_6 = \varphi(IpP_6) = \varphi([15, 26])$ = $\{ 17, 19, 23\}$, $MpP_6 = 23$, 

\hspace{.6cm} $p_7 = 17$, $p_8 = 19$, $p_9 = 23$.

$\Rightarrow$ $IpP_7 = [MIpP_6+1, 2MpP_6] = [27, 46]$, $MIpP_7 = 46$, 

\hspace{.5cm} $pP_7 = \varphi(IpP_7)$ $= \varphi([27, 46])$ $= \{ 29, 31, 37, 41, 43\}$, 

\hspace{.5cm}  $p_{10} = 29$, $p_{11} = 31$, $p_{12} = 37$, $p_{13} = 41$, $p_{14} = 43$, $MpP_7 = 43$.

$\Rightarrow$ $IpP_8 = [MIpP_7+1, 2MpP_7] = [47, 86]$, $MIpP_8 = 86$, 

\hspace{.5cm} $pP_8 = \varphi(IpP_8) = \varphi([47, 86])$ $= \{ 47, 53, 59, 61, 67, 71, 73, 79, 83\}$, 

\hspace{.5cm} $p_{15} = 47$, $p_{16} = 53$, $p_{17} = 59$, $p_{18} = 61$, $p_{19} = 67$, 

\hspace{.5cm} $p_{20} = 71$, $p_{21} = 73$, $p_{22} = 79$, $p_{23} = 83$, $MpP_8 = 83$. 

$\Rightarrow$ $IpP_9 = [MIpP_8+1, 2MpP_8] = [87, 166]$, $MIpP_9 = 166$, 

\hspace{.5cm} $pP_9 = \varphi(IpP_9) = \varphi([87, 166])$  $= \{89, 97, 101, 103, 107, 109,$ 

\hspace{.5cm} $113, 127, 131, 137, 139, 149, 151, 157, 163\}$, $MpP_9 = 163$, 

\hspace{.5cm} $p_{24} = 89$, $p_{25} = 97$, $p_{26} = 101$, $p_{27} = 103$, $p_{28} = 107$, 

\hspace{.5cm} $p_{29} = 109$, $p_{30} = 113$, $p_{31} = 127$, $p_{32} = 131$, $p_{33} = 137$, 

\hspace{.5cm} $p_{34} = 139$, $p_{35} = 149$, $p_{36} = 151$, $p_{37} = 157$, $p_{38} = 163$. 

$\Rightarrow$ $IpP_{10} = [MIpP_9+1, 2MpP_9] = [167, 326]$, $MIpP_{10} = 326$,

\hspace{.5cm} $pP_{10}$ $= \varphi(IpP_{10})$ $= \varphi([167, 326]) = \{167 = p_{39}, 173, 179, 181$,

\hspace{.5cm} $191, 193, 197, 199, 211, 223, 227, 229, 233, 239, 241, 251, 257, 263$,

\hspace{.5cm}  $269, 271, 277, 281, 283, 293, 307, 311, 313, 317 = p_{66}\}$, $\dots$.
\end{exm}

 With the notation, $pP_1 = \{2\}$, $pP_2 = \{3\}$, $pP_3 = \{5\}$, $pP_4 = \{7\}$, $pP_5 =$ $\{11, 13\}$, $pP_6 = \{17,19,23\}$, $\dots$, $IpP_{i+1} = [MIpP_{i}+1,$ $2MpP_{i}]$ and $pP_{i+1} = \varphi(IpP_{i+1})$ where $MIpP_{i} =$ max $IpP_{i} =$ the maximum value in $IpP_{i}$ and $MpP_i =$ $max~ pP_i$, $i = 2, 3, \ldots$, we obtain the following result. Hereafter, we call $pP_{i}$ as \textit{the $i^{th}$ pocket of prime(s)}, $i \in \mathbb{N}$.   

\begin{theorem} \label{c3} {\rm Let $pP_1 = \{2\}$, $IpP_1 = \{2\}$, $IpP_2 = [3, 4]$, $pP_2 = \{3\}$, $IpP_3 = [5, 6]$, $pP_3 = \{5\}$, $IpP_4 = [7, 10]$, $pP_4 = \{7\}$, $IpP_5 = [11, 14]$, $pP_5 = \{11, 13\}$, $\dots$ where $IpP_{i+1} = [MIpP_{i}+1,$ $2MpP_{i}]$, $pP_{i+1} = \varphi(IpP_{i+1})$, $MIpP_{i} =$ max $IpP_{i} =$ the maximum value in $IpP_{i}$ and $MpP_{i} =$ max $pP_{i} =$ the maximum value in $pP_{i}$, $i = 2, 3, \ldots$.~Then, 
\begin{enumerate} 
\item [\rm (a)]  The set of all prime numbers, 

$\mathbb{P}$ = $\varphi(\mathbb{N}) =$ $pP_1 \cup pP_2 \cup pP_3 \cup \ldots$ $=$ $\bigcup^{\infty}_{j=1}$ $pP_j$.

\item [\rm (b)]  The set of all pockets of primes, 

$\{pP_1, pP_2, \ldots, pP_{i}, \ldots \}$  partitions the set of all primes. \hfill $\Box$
\end{enumerate} }
\end{theorem}

\begin{obs}\quad  {\rm
\begin{enumerate} 
\item [\rm (i)]  To obtain the set of all primes $\mathbb{P}$, one can consider different sets of pockets of primes. 
\item [\rm (ii)] A set of pockets of primes whose union is the set of all primes need not be a partition of $\mathbb{P}$.  
\end{enumerate} }
\end{obs} 
We generated prime numbers from successive intervals $IpP_{j+1}$: from each successive interval, we calculate all prime numbers and thereby obtain corresponding pocket of primes $pP_{j+1}$ = $\varphi([MIpP_j+1, 2MpP_j])$ for $j = 1, 2, \ldots$, starting with $pP_1$ = $\{2\}$, $IpP_1 = \{2\}$ and then find $IpP_{j+1}$ = $[MIpP_{i}+1,$ $2MpP_{i}]$, $MIpP_{j}$ = max $IpP_{j}$ and $MpP_{j}$ = max $pP_{j}$. Now, the question is `Is it possible to find a better method to generate primes?'. The following method based on $IpP_{j+1}$ = $[MIpP_{j}+1, MpP_j^2]$ and $pP_{j+1}$ = $\varphi(IpP_{j+1})$ = $\varphi([MIpP_{j}+2, MpP_j^2-2])$ where $MIpP_{j}$ = maximum value in $IpP_{j}$ and $MpP_{j}$ = maximum value in $pP_{j}$ with $p_1$ = 2, $pP_1$ = $\{2\}$, $IpP_2$ = $\{3\}$ = $pP_2$ and $p_2$ = 3 = $MIpP_2$ = $MpP_2$ for $j = 2, 3, \ldots$ is a better method. To begin with, let us calculate prime numbers in $[p_j+1, p_j^2],$ $j\geq 2$.   

{\subsection ~  {\bf Calculation of Prime numbers in $[p_{j}+1, p_j^2],$~ $j\geq 2$}}

Here, we obtain prime numbers by considering successive intervals $IpP_{j+1}$ as $[p_{j}+1, p_j^2]$~ instead of $[p_j+1, 2p_j]$ in the  previous method for $j\geq 2$ with $p_1 = 2$, $p_2 = 3$, $\ldots$. Interval $[p_j+1, p_j^2]$, in general, contains more primes than in $[p_j+1, 2p_j]$ for $j\geq 2$. Similar to Theorems \ref{c2} and \ref{c3}, we obtain Theorems \ref{c4} and \ref{c5} corresponding to intervals $[p_j+1, p_j^2]$, $j \geq 2$. 

\begin{theorem} \label{c4} {\rm Let $1 \leq i \leq j$, $i,j\in\mathbb{N}$, $Q_{i,j} = \left\lfloor \frac{p_j}{p_i} \right\rfloor$ and $Q_{i,j}^{''} = \left\lfloor \frac{p_j^2-2}{p_i} \right\rfloor$ and $S_j$ =  $[p_j+1, p_j^2] \setminus \mathbb{C}([p_j+1, p_j^2])$. Then for $j \geq 2$, 
\begin{enumerate} 
\item [\rm (a)] \label{c41} $S_j$ is non-empty; 
\item [\rm (b)] \label{c42} $S_j$ contains prime number(s) as its element(s); 
\item [\rm (c)] \label{c43} $S_j$ = $\varphi([p_j+1, p_j^2])$ = $\varphi([p_j+2, p_j^2-2])$; 
\item [\rm (d)] \label{c44} $S_j$ = $[p_j+2, p_j^2-2] \setminus (\bigcup^{j}_{i=1} p_i\{Q_{i,j}+1,Q_{i,j}+2,\ldots,Q_{i,j}^{''}\})$.
\end{enumerate} }
\end{theorem} 
\begin{proof} By  Bertrand postulate, $[p_j+1, 2p_j]$ contains at least one prime, $j\in\mathbb{N}$. For $j\in \mathbb{N}$, $[p_j+1, p_j^2]$ $\supseteq$ $[p_j+1, 2p_j]$ and so $[p_j+1, p_j^2]$ contains at least one prime since $p_j \geq 2$. Also, for $j \geq 2$, $p_j$ is odd prime and $p_j^2$ an odd composite number. And hence the set of all prime numbers contained in $[p_j+1, p_j^2]$ is same as the set of all prime numbers contained in $[p_j+2, p_j^2-2]$ and is obtained by removing all composite numbers contained in $[p_j+2, p_j^2-2]$. Hence $(a)$, $(b)$ and $(c)$ are true.

Using $(a)$, $[p_j+2, p_j^2-2]$ contains at least one prime for $j \geq 2$. And any prime number contained in $[p_j+1, p_j^2]$ is greater than $p_j$ and less than $p_j^2$ for $j \geq 2$. Also, for $k,l\in \mathbb{N}$, prime numbers $p_{j+k},p_{j+l} > p_j$ and $p_{j+k}p_{j+l} > p_j^2$. This implies, $p_1$, $p_2$, $\ldots$, $p_{j-1}$ or $p_j$ is a divisor of every composite number contained in $[p_j+1, p_j^2]$ and after removal of all  composite numbers that are multiples of $p_1$, $p_2$, $\ldots$, $p_{j-1}$ or $p_j$ from $[p_j+1, p_j^2]$, the resultant set contains prime(s) only. This implies,  

$[p_j+1, p_j^2]$ $\setminus$ $\mathbb{C}([p_j+1, p_j^2])$ = $\varphi([p_j+1, p_j^2])$ = $\varphi([p_j+2, p_j^2-2])$

= $[p_j+2, p_j^2-2]$ $\setminus$ $\mathbb{C}([p_j+2, p_j^2-2])$ since $p_j+1$ is even for $j\geq 2$.

= $[p_j+2, p_j^2-2]$ $\setminus$ ($\bigcup^{j}_{i=1}$ $p_i\{Q_{i,j}+1,Q_{i,j}+2,\ldots,Q_{i,j}^{''}\})$ by Lemma \ref{b1}(f) where $Q_{i,j} = \left\lfloor \frac{p_j}{p_i} \right\rfloor$ and $Q_{i,j}^{''} = \left\lfloor \frac{p_j^2-2}{p_i} \right\rfloor$, $1 \leq i \leq j$.

Hence, result $(d)$ is true.
\end{proof}

In the above result, $\varphi([p_j+1, p_j^2])$, the set of all prime number(s) in the interval $[p_j+1, p_j^2]$, is a non-empty set for every $j$ using Bertrand postulate, $j\in \mathbb{N}$. Now, let us calculate different pockets of primes $pP_{j+1}$ = $\varphi([MIpP_j+2, MpP_j^2-2])$ with $p_1$ = 2, $IpP_1$ = $\{2\}$ = $pP_1$, $\pi(pP_1)$ = $|pP_1|$ = 1, $MIpP_1$ = 2 = $MpP_1$, $IpP_2$ = $[MIpP_1+1$, $MpP_1^2]$ = $[3, 4]$, $pP_2$ = $\varphi(IpP_2)$ = $\{3 \}$,  $\pi(pP_2)$ = 1, $MIpP_2$ = 4, $MpP_2$ = 3 = $p_{\pi(pP_1)+\pi(pP_2)}$ = $p_2$, $IpP_{j+1}$ = $[MIpP_j+1, MpP_j^2]$, $pP_{j+1}$ = $\varphi([MIpP_j+2, MpP_j^2-2])$, $\pi(pP_{j+1})$ = $|pP_{j+1}|$ and $MpP_{j+1}$ = max $pP_{j+1}$ = $p_{\pi(pP_1)+\pi(pP_2)+\ldots+\pi(pP_j)+\pi(pP_{j+1})}$, $j \geq 2$ and $j\in \mathbb{N}$. Here, we obtain $pP_j$ from $IpP_j$ using Theorem \ref{c4}$(d)$.

\begin{exm} Let $p_1 = 2$, $IpP_1 = \{2\} = pP_1$, $\pi(pP_1) = 1 = |pP_1|$, $MIpP_1 = 2 = MpP_1$.

$\Rightarrow$ $IpP_2 = [MIpP_1+1, MpP_1^2] = [3, 4]$, $pP_2 = \varphi(IpP_2) = \{3 \}$, 
 
\hspace{.4cm} $\pi(pP_2) = |pP_2| = 1$, $MpP_2 = 3 = p_{\pi(pP_1)+\pi(pP_2)} = p_2$, $MIpP_2 = 4$.

$\Rightarrow$ $IpP_3 = [MIpP_2+1, MpP_2^2] = [5, 9]$, $MIpP_3 = 9$, 
 
\hspace{.5cm} $pP_3 = \varphi(IpP_3) = \{5, 7\}$, $\pi(pP_3) = |pP_3| = 2$, $p_3$ = 5, $p_4$ = 7,   

\hspace{.5cm} $MpP_3 = 7 = p_{\pi(pP_1)+\pi(pP_2)+\pi(pP_3)} = p_4$.

$\Rightarrow$  $IpP_4 = [MIpP_3+1, MpP_3^2]$ $= [10, 49]$, $MIpP_4 = 49$,  

\hspace{.5cm} $pP_4 = \varphi(IpP_4)$  = $\varphi([10, 49])$ = $\varphi([11, 47])$ 

\hspace{1.2cm} = $\{ 11, 13, 17, 19, 23, 29, 31, 37, 41, 43, 47\}$, $\pi(pP_4)$ = $|pP_4|$ = $11$,

\hspace{.5cm} $MpP_4$ = 47 = $p_{\pi(pP_1)+\pi(pP_2)+\pi(pP_3)+11}$ = $p_{4+11} = p_{15}$, 

\hspace{.5cm} $p_5$ = 11, $p_6$ = 13, $p_7 = 17$, $p_8 = 19$, . . . ,  $p_{14} = 43$, $p_{15}$ = 47.    

$\Rightarrow$  $IpP_5 = [MIpP_4+1, MpP_4^2]$ $= [50, 2209]$, $MIpP_5 = 2209$,  

\hspace{.5cm} $pP_5 = \varphi(IpP_5) = \varphi([MIpP_4+1, MpP_4^2])$ $= \varphi([51, 2209])$   

\hspace{1.3cm}  $= \varphi([53, 2207])$ = $\{ 53, 59, 61, 67, 71, 73, 79, 83, 89, 97$,  

$101, 103, 107, 109, 113, 127, 131, 137, 139, 149, 151, 157, 163, 167, 173$, 

$179, 181, 191, 193, 197, 199, 211, 223, 227, 229, 233, 239, 241, 251, 257$, 

$263, 269, 271, 277, 281, 283, 293, 307, 311, 313, 317, 331, 337, 347, 349$, 

$353, 359, 367, 373, 379, 383, 389, 397, 401, 409, 419, 421, 431, 433, 439$, 

$443, 449, 457, 461, 463, 467, 479, 487, 491, 499, 503, 509, 521, 523, 541$, 

$547, 557, 563, 569, 571, 577, 587, 593, 599, 601, 607, 613, 617, 619, 631$, 

$641, 643, 647, 653, 659, 661, 673, 677, 683, 691, 701, 709, 719, 727, 733$, 

$739, 743, 751, 757, 761, 769, 773, 787, 797, 809, 811, 821, 823, 827, 829$, 

$839, 853, 857, 859, 863, 877, 881, 883, 887, 907, 911, 919, 929, 937, 941, 947$, 

$953, 967, 971, 977, 983, 991, 997, 1009, 1013, 1019, 1021, 1031, 1033, 1039$, 

$1049, 1051, 1061, 1063, 1069, 1087, 1091, 1093, 1097, 1103, 1109, 1117$, 

$1123, 1129, 1151, 1153, 1163, 1171, 1181, 1187, 1193, 1201, 1213, 1217$, 

$1223, 1229, 1231, 1237, 1249, 1259, 1277, 1279, 1283, 1289, 1291, 1297$, 

$1301, 1303, 1307, 1319, 1321, 1327, 1361, 1367, 1373, 1381, 1399, 1409$, 

$1423, 1427, 1429, 1433, 1439, 1447, 1451, 1453, 1459, 1471, 1481, 1483$, 

$1487, 1489, 1493, 1499, 1511, 1523, 1531, 1543, 1549, 1553, 1559, 1567$, 

$1571, 1579, 1583, 1597, 1601, 1607, 1609, 1613, 1619, 1621, 1627, 1637$, 

$1657, 1663, 1667, 1669, 1693, 1697, 1699, 1709, 1721, 1723, 1733, 1741$, 

$1747, 1753, 1759, 1777, 1783, 1787, 1789, 1801, 1811, 1823, 1831, 1847$, 

$1861, 1867, 1871, 1873, 1877, 1879, 1889, 1901, 1907, 1913, 1931, 1933$, 

$1949, 1951, 1973, 1979, 1987, 1993, 1997, 1999, 2003, 2011, 2017, 2027$, 

$2029, 2039, 2053, 2063, 2069, 2081, 2083, 2087, 2089, 2099, 2111, 2113$, 

$2129, 2131, 2137, 2141, 2143, 2153, 2161, 2179, 2203, 2207\}$, $\pi(pP_5) = 314$, 

\hspace{.2cm} $p_{15+1} = p_{16} = 53$, $p_{17} = 59$, $\ldots$, $p_{15+\pi(pP_5)} = p_{329} = MpP_5 = 2207$. 

$\Rightarrow$ $IpP_{6} = [MIpP_5+1, MpP_5^2]$ $= [2210, 4870849])$,  $MIpP_6 = 4870849$,

\hspace{.5cm} $pP_{6} = \varphi([MIpP_5+2, MpP_5^2-2])$ $= \varphi([2211, 4870847])$

\hspace{1.5cm} $= \{2213 = p_{330}, 2221 = p_{331}, \ldots, 4870843 = p_{340059}\}$,   

\hspace{.5cm} $\pi(pP_6) = 339730$, $MpP_6$ = 4870843 = $p_{329+|pP_6|} = p_{340059}$, $\dots$.
\end{exm}

Continuing the above process of obtaining prime numbers from pockets of primes, one can generate consecutive prime numbers up to any given limit. The algorithm used here is $p_1$ = 2, $IpP_1$ = $\{2\}$ = $pP_1$, $\pi(pP_1)$ = $|pP_1|$ = 1, $MIpP_1$ = 2 = $MpP_1$, $p_2$ = 3, $IpP_2$ = $\{3\}$ = $pP_2$, $MIpP_2$ = 3 = $MpP_2$, $\pi(pP_2)$ = 1, then $IpP_{i+1}$ = $[MIpP_{i}+1, MpP_{i}^2]$, $pP_{i+1}$ = the ${(i+1)}^{th}$ pocket of primes = $\varphi([MIpP_{i}+1, MpP_{i}^2])$ = $\varphi([MIpP_{i}+2, MpP_{i}^2-2])$, $\pi(pP_{i+1})$ = $|pP_{i+1}|, MpP_{i}$ = max $pP_{i}$ = the maximum value of $pP_{i}$ = the biggest prime in $pP_{i}$ = the biggest prime in $IpP_{i}$, $i$ = $2, 3, \ldots$.  
  
\begin{rem} The above algorithm works since $MpP_i^2$  is composite and odd for every $i \geq 2$ whereas $MpP_1^2 = 4$ is even. Thus, corresponding to the above algorithm, we get the following result.  
\end{rem}

\begin{theorem} \label{c5} \quad {\rm Let $p_1$ = 2, $IpP_1$ = $\{2\}$ = $pP_1$, $p_2$ = 3, $IpP_2$ = $\{3\}$ = $pP_2$, $p_2$ = 3 = $MIpP_2$ = $MpP_2$, $IpP_{i+1}$ = $[MIpP_{i}+1, MpP_{i}^2]$ and $pP_{i+1}$ = $\varphi([MIpP_{i}+1, MpP_{i}^2])$ = $\varphi([MIpP_{i}+2, MpP_{i}^2-2])$ where $MIpP_{i}$ = max $IpP_{i}$ and $MpP_{i}$ = max $pP_{i}$, $i$ = $2, 3, \ldots$. Then, the pockets of primes are 
\begin{enumerate}
\item [\rm (a)] $pP_1$ = $\varphi(\{2\})$ = $\{2\}$, $pP_2$ = $\varphi(\{3\})$ = $\{3\}$, $pP_3$ = $\{5, 7\}$,

\noindent
 $pP_4$ = $\{ 11 = p_5, 13, 17, 19, 23, 29, 31, 37, 41, 43, 47 = p_{15}\}$;

\noindent
 $pP_5$ = $\{53 = p_{16}, 59 = p_{17}, 61, \ldots, 2207 = p_{329} = MpP_5\}$, 

\noindent
 $pP_{6} = \{2213 = p_{330}, 2221 = p_{331}, \dots, 4870843 = p_{340059}\}$, 

$\ldots$. 
\item  [\rm (b)]  The set of all prime numbers, 

$\mathbb{P}$ = $\varphi(\mathbb{N})$ = $\bigcup_{j=1}^\infty pP_{j}$ = $pP_1 \bigcup pP_2$ $\bigcup$  $(\bigcup_{j=2}^\infty \varphi([MIpP_j+1, MpP_j^2]))$

\hspace{.25cm} = $pP_1 \bigcup pP_2$ $\bigcup$  $(\bigcup_{j=2}^\infty \varphi([MIpP_j+2, MpP_j^2-2]))$ 

\hspace{.25cm} = $pP_1 \bigcup pP_2$ $\bigcup (\lim_{n \rightarrow \infty} (\bigcup_{j=2}^n \varphi([MIpP_j+2, MpP_j^2-2])))$.
\item  [\rm (c)] The set of all pockets of primes, $\{pP_1, pP_2, \ldots\}$ partitions the set of all prime numbers, $\mathbb{P}$. 
\end{enumerate} }
\end{theorem}
\begin{proof} Here, for $i \neq j$ and $i,j\in \mathbb{N}$, $pP_i \cap pP_j$ = $\emptyset$ and each pocket of prime(s) is non-empty by Bertrand Postulate. Also, the set of all intervals $IpP_j$ partitions $\mathbb{N}\setminus \{1\}$ and each interval $IpP_j$ covers pocket of prime(s) $pP_j$, $j$ = $1, 2, \ldots$. Then, the theorem is true from the following.

$\mathbb{P}$ = $\varphi(\mathbb{N})$ = $\bigcup_{j=1}^\infty pP_{j}$ = $\lim_{n \rightarrow \infty} (\bigcup_{j=1}^n pP_j)$ = $pP_1 \bigcup pP_2$ $\bigcup (\bigcup_{j=3}^n pP_j)$ 

\hspace{.25cm} = $pP_1 \bigcup pP_2 \bigcup$ $(\lim_{n \rightarrow \infty} (\bigcup_{j=2}^n \varphi([MpP_j+1, MpP_j^2])))$.

\hspace{.25cm} = $pP_1 \bigcup pP_2 \bigcup$ $(\lim_{n \rightarrow \infty} (\bigcup_{j=2}^n \varphi([MpP_j+2, MpP_j^2-2])))$.
\end{proof}

{\subsection ~  {\bf Calculation of Prime numbers in $[p_{j}+1, p_j^2+4p_j+3]$, $j\geq 2$}}

 In the two methods that we have discussed to obtain prime numbers using successive pockets of prime(s) $pP_1$, $pP_2$, $\ldots$ with $IpP_1$ $= \{2\}$ $= pP_1$, $IpP_2 =$ $\{3\}$ $= pP_2$, $MIpP_j = max ~ IpP_j$ and $MpP_j$ $= max ~ pP_j$, in the first method we take $IpP_{j+1} =$ $[MIpP_j+1,$ $2MpP_j]$ and $pP_{j+1} = \varphi(IpP_{j+1}) = \varphi([MIpP_j$ $+1,$ $2MpP_j])$ for $j = 2,3,\ldots$ and in the second method we take $IpP_{j+1} =$ $[MIpP_j+1, MpP_j^2]$ and $pP_{j+1} = \varphi(IpP_{j+1}) =$ $\varphi([MIpP_j+2,$ $MpP_j^2-2])$ for $j = 2,3,\ldots$. Let $\pi(pP_j) =|pP_j|$, $j = 1, 2, \ldots$. It is easy to observe the following.

\begin{enumerate}
\item [\rm (i)] The second method is better than the first in terms of number of primes generated in successive pockets of primes. The number of elements in the successive pockets of primes in the two methods are as follows. 

\noindent
{\bf Method-1:} $\pi(pP_1) = 1$, $\pi(pP_2) = 1$, $\pi(pP_3) = 1$, $\pi(pP_4) = 1$, $\pi(pP_5) = 2$, $\pi(pP_6) = 3$, $\pi(pP_7) = 5$, $\pi(pP_8) = 9$, $\pi(pP_9) = 15$, $\pi(pP_{10}) = 28$, $\ldots$.  

\noindent
{\bf Method-2:} $\pi(pP_1) = 1$, $\pi(pP_2) = 1$, $\pi(pP_3) = 2$, $\pi(pP_4) = 11$, $\pi(pP_5) = 314$, $\pi(pP_6) = 339730$, $\ldots$. 

See Table $1$ for more details.
\item [\rm (ii)]  To generate prime numbers, if we consider any bigger interval of the form $[MIpP_j+1,$ $MpP_j^2+k]$ in the previous methods, then the method may fail since $[MIpP_j+1, MpP_j^2+k]$ may contain composite numbers of the form $MpP_j^2+m$ for some $k$ and $m$ such that $1 \leq m \leq k$ and each of its prime divisor is $> MpP_j$. 

\item [\rm (iii)]  In the previous remark, we raised the question whether there exists such a value of $k$ for which the algorithm works in general to generate prime numbers? If it exists, is it possible to find out such value(s) of $k$? And what is the maximum general value of $k$ ? Yes, such values of $k$ exist. See the following. 

\item [\rm (iv)]  For $p_1 = 2$, $p_2 = 3$, $pP_1 = \{2, 3\}$, $IpP_1 = [2, 3]$ and $j \geq 2$, $p_j$ is an odd prime and the next prime $p_{j+1} \geq p_j+2$ and so $p_{j+1}p_{j+1}$ $\geq$ ${(p_j+2)}^2$ $= p_j^2 + 4p_j +4$ $> p_j^2 + 4p_j +3$. This implies, $k$ $= p_j^2+4p_j+3$ is a possible value of $k$, in general, for which the following algorithm works.  

$p_1 = 2$, $p_2 = 3$, $IpP_1 = [2, 3]$, $pP_1 = \{2, 3\}$,

$\pi(pP_1) = 2$, $MIpP_1 = 3 = MpP_1 = p_2$, 

$IpP_{j+1} = [MIpP_j+1, MpP_j^2+k]$, 

$pP_{j+1} = \varphi([MIpP_j+1, MpP_j^2+k-1])$, 

$\pi(pP_{j+1}) =$ $|pP_{j+1}|$, $MIpP_{j+1} =$ $max~ IpP_{j+1}$ and 

$MpP_{j+1} =$ $max~pP_{j+1}$ $= p_{\pi(pP_1)+\pi(pP_2)+\ldots+\pi(pP_{j+1})}$ for $j = 1, 2, \ldots$ 

where $k = 4MpP_j+3$ and $MpP_j^2+4MpP_j+3$ is even. 
\item [\rm (v)]  Continuing the above process of obtaining prime numbers from pockets of primes, one can generate consecutive prime numbers up to any given limit. Thus, corresponding to the third method, we get the following important result. 
\end{enumerate}   

\begin{rem} In the above algorithm if we take $k = MpP_j^2+4MpP_j+4$ instead of $MpP_j^2+4MpP_j+3$, then the algorithm fails when $MpP_j$ and its successive prime are twin primes (two primes of difference two). Thus, $k = MpP_j^2+4MpP_j+3$ is a possible, in general, maximum value of $k$ used in the third method to generate successive primes and correspondingly we get the most important result, Theorem \ref{c6}. And to calculate successive pockets of primes $pP_j$ from $IpP_j$, we use $(c)$ of Theorem \ref{c7}.
\end{rem}  

\begin{theorem} \label{c6} {\rm  For $j\in \mathbb{N}$ and $k\in \mathbb{N}_0$, the following statements are true. 
\begin{enumerate}
\item [\rm (a)] $k$ = $4MpP_j+3$ is a possible, in general, maximum value of $k$ for which the following algorithm works to generate successive prime numbers from successive pockets of primes. 

$p_1$ = 2, $p_2$ = 3, $IpP_1$ = $[2, 3]$, $pP_1$ = $\{2, 3\}$, $\pi(pP_1)$ = 2, 

 $MpP_1$ = 3 = $MIpP_1$ = $p_{\pi(pP_1)}$ = $p_2$,  

$IpP_{j+1}$ = $[MIpP_{j}+1, MpP_{j}^2+k]$, $MIpP_j$ = max $IpP_j$, 

$pP_{j+1}$ = $\varphi(IpP_{j+1})$ = $\varphi([MIpP_{j}+1, MpP_{j}^2+k-1])$, 

$\pi(pP_{j+1})$ = $|pP_{j+1}|$ and 

$MpP_{j+1}$ = max  $pP_{j+1}$ = $p_{\pi(pP_1)+\pi(pP_2)+\ldots+\pi(pP_{j+1})}$ for $j$ = $1, 2, \ldots$.  

\item [\rm (b)] Let $IpP_1$ = $[2, 3]$, $pP_1$ = $\{2, 3\}$, $p_1$ = 2, $p_2$ = 3, $\pi(pP_1)$ = 2 and $MIpP_1$ = 3 = $MpP_1$ = $p_{\pi(pP_1)}$ = $p_2$. Then, the set of all primes, 
\\
 $\mathbb{P}$ = $\bigcup_{j=1}^{\infty} pP_{j}$
 = $\lim_{n \rightarrow \infty} (\bigcup_{j=1}^n pP_j)$

= $pP_1$ $\bigcup (\lim_{n\rightarrow\infty} (\bigcup_{j=1}^n \varphi([MIpP_j+1, MpP_j^2+4MpP_j+2])))$ 
\\
where $IpP_{j+1}$ = $[MIpP_j+1, MpP_j^2+4MpP_j+3]$, $MIpP_{j}$ = max $IpP_j$, $pP_{j+1}$ = $\varphi(IpP_{j+1})$ = $\varphi([MIpP_{j}+1, MpP_j^2+4MpP_j+2])$, $\pi(pP_{j+1})$ = $|pP_{j+1}|$ and $MpP_{j+1}$ = $p_{\pi(pP_1)+\pi(pP_2)+\ldots+\pi(pP_{j+1})}$ = max $pP_{j+1}$, $j$ = $1, 2, \ldots$. \hfill $\Box$
\end{enumerate} }
\end{theorem}

\begin{definition} Let $pP_1 = \{2, 3\}$, $p_1 = 2$, $p_2 = 3$, $\pi(pP_1) = 2$, $IpP_1 = [2, 3]$, $MIpP_1$ $= 3$ $= MpP_1$, $IpP_{j+1} = [MIpP_j+1,$ $MpP_j^2$ $+4MpP_j+3]$, $MIpP_{j+1} =$ $max ~ IpP_{j+1}$, $pP_{j+1} = \varphi([MIpP_j+1, MpP_j^2$ $+4MpP_j+2])$, $\pi(pP_{j+1}) =$ $|pP_{j+1}|$ and $MpP_{j+1} =$ $p_{\pi(pP_1)+\pi(pP_2)+\ldots+\pi(pP_{j+1})} =$ $max$ $pP_{j+1}$, $j = 1, 2, \ldots$. Then, $pP_1$, $pP_2$, $\ldots$ are called {\em the standard pockets of primes}.   
\end{definition}
\begin{theorem} \label{c7} {\rm Let $j \geq 2$ and $S_j$ = $[p_j+1, p_j^2+4p_j+3] \setminus \mathbb{C}([p_j+1, p_j^2+4p_j+3])$. Then  
\begin{enumerate} 
\item [\rm (a)] \label{c71} $S_j$ is non-empty; 
\item [\rm (b)] \label{c72} $S_j$ contains prime number(s) as its element(s); 
\item [\rm (c)] \label{c73} $S_j$ = $\varphi([p_j+1, p_j^2+4p_j+3])$ =  $\varphi([p_j+2, p_j^2+4p_j+2])$; 
\item [\rm (d)] \label{c74} $S_j$ = $[p_j+2, p_j^2+4p_j+2] \setminus$ $(\bigcup^{j}_{i=1} p_i\{Q_{i,j}+1,Q_{i,j}+2,\ldots,Q_{i,j}^{'''}\})$ 

where $Q_{i,j}$ = $\left\lfloor \frac{p_j}{p_i} \right\rfloor$ and $Q_{i,j}^{'''}$ = $\left\lfloor \frac{p_j^2+4p_j+2}{p_i} \right\rfloor$, $1 \leq i \leq j$.
\end{enumerate} }
\end{theorem} 
\begin{proof} Proof is similar to Theorem \ref{c4}.
\end{proof}

Now, we calculate the standard pockets of primes using Theorems \ref{c6} and \ref{c7}.
\begin{exm} $p_1 = 2$, $p_2 = 3$, $IpP_1 = [2, 3]$, $pP_1 = \{2, 3\}$, $MIpP_1 = 3 = MpP_1$, $\pi(pP_1) =2$;

$\Rightarrow IpP_2 = [MIpP_1+1, MpP_1^2+4MpP_1+3] = [4, 24]$, $MIpP_2 = 24$,

\hspace{.5cm} $pP_2 = \varphi([MIpP_1+2, MpP_1^2+4MpP_1+2])$ $= \varphi([5, 23])$ $= \{5, 7$,

\hspace{.5cm}  $11, 13, 17, 19, 23\}$, $\pi(pP_2) = 7$, $p_3 = 5$, $p_4 = 7$, $p_5 = 11$, $p_6 = 13$,  

 \hspace{.5cm}  $p_7 = 17$, $p_8 = 19$, $p_9 = 23$, $MpP_2 = 23 = p_{\pi(pP_1)+\pi(pP_2)} = p_{9}$; 

 $\Rightarrow$ $IpP_3 = [MIpP_2+1, MpP_2^2+4MpP_2+3] = [25, 624]$, 

\hspace{.1cm} $MIpP_3 =$ $624$, $pP_3 = \varphi([MIpP_2+2, MpP_2^2+4MpP_2+2])$ 

\hspace{.1cm} $= \varphi([26, 623]) = \{29, 31, 37, 41, 43, 47, 53, 59, 61, 67, 71, 73, 79, 83$,     

$89, 97, 101, 103, 107, 109, 113, 127, 131, 137, 139, 149, 151, 157, 163$, 

$167, 173, 179, 181, 191, 193, 197, 199, 211, 223, 227, 229, 233, 239, 241$, 

$251, 257, 263, 269, 271, 277, 281, 283, 293, 307, 311, 313, 317, 331, 337$, 

$347, 349, 353, 359, 367, 373, 379, 383, 389, 397, 401, 409, 419, 421, 431$, 

$433, 439, 443, 449, 457, 461, 463, 467, 479, 487, 491, 499, 503, 509, 521$, 

$523, 541, 547, 557, 563, 569, 571, 577, 587, 593, 599, 601,  607, 613, 617$, 

$619\}$, ~~ $\pi(pP_3) = 105$,

$MpP_3 = 619 = p_{9+\pi(pP_3)} = p_{114}$, $p_{10} = 29$, $p_{11} = 31$, $\ldots$, $p_{114} = 619$;

$\Rightarrow$ $IpP_4 = [MIpP_3+1, MpP_3^2+4MpP_3+3] = [625, 385640]$, 

\hspace{.5cm} $MIpP_4 =$ $385640$, $pP_4 = \varphi([MIpP_3+2, MpP_3^2+4MpP_3+2])$

\hspace{.5cm}  $= \varphi([626, 385639]) = \{631, 641, \ldots, 385639\}$, 

\hspace{.5cm} $\pi(pP_4) = 32622$, $MpP_4 =$ $385639 = p_{114+|pP_4|} = p_{32736}$, 

\hspace{.5cm} $p_{115} = 631$, $p_{116} = 641$, $\ldots$, $p_{32736} = 385639 = MpP_4$;
 
 $\Rightarrow$ $IpP_5 = [MIpP_4+1, MpP_4^2+4MpP_4+3] $ = $[385641, 148718980880]$, 

\hspace{1.2cm}  $MIpP_5 =$ $148718980880$, 

\hspace{.5cm} $pP_5 = \varphi([MIpP_4+2, MpP_4^2+4MpP_4+2])$ 

\hspace{1.2cm} $= \varphi([385642,$ $148718980879])$ $= \{385657 = p_{32737}$, 

\hspace{1.2cm} $385661 = p_{32738}, \ldots, MpP_5 \leq 148718980879\}$,   

\hspace{.5cm} $\pi(pP_5) = |pP_5|$, $MpP_5 =$ $p_{32736+\pi(pP_5)}$; ~~ $\ldots$.
\end{exm}

\section{Algorithm to generate primes}
To generate prime numbers to any large extend, the third method is the fastest and the best among the three methods that we have discussed. Here, we present algorithm corresponding to the third method only that is based on Lemma \ref{b1} and Theorems \ref{c6} and \ref{c7}. It is easy to check whether a given natural number is prime or not by comparing with the list of primes already generated provided the number is less than or equal to any prime that was already generated. 

\vspace{.3cm}
\noindent
{\bf Generating primes with $IpP_{j+1}$ = $[MIpP_j+1, MpP_j(MpP_j+4)+3]$
\\
Algorithm:}

$p_1 = 2$, $p_2 = 3$, $IpP_1 = [2, 3]$, $pP_1 = \{2, 3\}$,  

$\pi(pP_1) = 2$, $MIpP_1 = 3 = MpP_1$, 

$j$ $=$ $2$, ~~~(Number of primes already generated)

$x = 8$, ~(No. of pockets of primes we generate in this particular program.)

Do $10$ $k = 1$ to $x$,

$IpP_{k+1} = [MIpP_k+1, MpP_k(MpP_k+4)+3]$,

$IP_{k+1} = [MIpP_k+2, MpP_k(MpP_k+4)+2]$,

$C(IP_{k+1}) = \{\}$,

Do $20$ $i = 1$ to $j$, 

$Q_{i,k} = \left\lfloor \frac{MpP_k}{p_i}\right\rfloor$, 

$Q_{i,k}^{'''} = \left\lfloor \frac{MpP_k^2+4MpP_k+2}{p_i}\right\rfloor$,

$C_i(IP_{k+1}) = \{p_i(Q_{i,k} +1), p_i(Q_{i,k} +2), \ldots, p_iQ_{i,k}^{'''}\}$,

$20$ $C(IP_{k+1})$ $=$ $C(IP_{k+1})$ $\bigcup$ $C_i(IP_{k+1})$

$pP_{k+1}$ $=$  $IP_{k+1}$ $\setminus$ $C(IP_{k+1})$,~~($C(IP_{k})$ contains all composite No.s of $IP_k$.)

$\pi(pP_{k+1})$ $=$ $|pP_{k+1}|$, 

$MpP_{k+1} =$ $max~ pP_{k+1}$,

$tpP_{k+1}$ $=$ $pP_{k+1}$, ~~(Here, $tpP_{k+1}$ represents temporary $pP_{k+1}$.)

$y$ $=$ $j$,

$j$ $=$ $j+\pi(pP_{k+1})$,

$p_{j} =$ $MpP_{k+1}$,

print (',$k+1$, $^{th}$ pocket of primes with No. of primes =, ', $\pi(pP_{k+1})$, 

`starting with', $y+1$, $^{th}$ prime to ', $j$, $'^{th}$ prime. They ~are') 

print $(p_{y+l}, ~ l = 1$ to $\pi(pP_{k+1})$, 

$10$ Print $('MpP_{k+1} =', p_j, '=', j,$ $'^{th}$~ prime.')

End 

------------------------------------------

In the three methods that we have discussed we could see that the values of $\pi(pP_k)$s, in each method, play an important roll and thereby we define the following. 
\begin{definition} The sequence $\pi(pP_1)$, $\pi(pP_2)$, $\ldots$ is called {\textit the sequence of orders of pockets of primes} $pP_1$, $pP_2$, $\ldots$.
\end{definition}
Table $1$ shows different sequences of orders of pockets of primes in the three methods up to $k$ $=$ $8$ in the first method, $k = 6$ in the second and $k = 4$ in the third. 

\begin{rem} A new study is needed on the behaviour of sequences of orders of pockets of primes in the three methods. 
\end{rem}
\begin{rem} While generating prime numbers using the third method, one need not start with $IpP_1$ = $[2, 3]$, $pP_1$ = $\{2, 3\}$, $IpP_2$ = $[4, 24]$, $\ldots$.

If $p_1$ = 2, $p_2$ = 3, $\ldots$, $p_k$ are already known primes, then by taking 

$pP_1$ = $\{p_1$ = 2, $p_2$ = 3, $\ldots, p_k\}$, $IpP_1$ = $[2, p_k]$, 

$MIpP_1$ = $p_k$ = $MpP_1$, $\pi(pP_1)$ = $k$, 

For $j = 1, 2, \ldots$, 

$IpP_{j+1}$ = $[MIpP_j+1, MpP_j^2+4MpP_j+3]$, 

$MIpP_{j+1}$ = max $IpP_{j+1}$, 

$pP_{j+1}$ = $\varphi(IpP_{j+1})$ = $\varphi([MIpP_j+1, MpP_j^2+4MpP_j+2])$, 

$\pi(pP_{j+1})$ = $|pP_{j+1}|$;

$MpP_{j+1}$ = max $pP_{j+1}$ = $p_{\pi(pP_1)+\pi(pP_2)+\ldots+\pi(pP_{j+1})}$.
\end{rem}

\section{Riemann zeta function and prime numbers generation}
In this section, using Riemann zeta function, we generate prime numbers. Riemann zeta function $\zeta(z)$ \cite{jps} is given by 

\hspace{.5cm} $\zeta(z) = \frac{1}{1^z}+\frac{1}{2^z}+\frac{1}{3^z}+\frac{1}{4^z}+\ldots$,
\\
the series is convergent for $Re(z) > 1$.

\begin{prop} \quad \label{d1} {\rm Let $p_1, p_2, \ldots$ denote the primes 2, 3, \ldots numbered in increasing order, $k\in\mathbb{N}$, $Re(z) > 1$,  $IpP_{k+1}$ = $[p_{k}+1, p_{k}(p_k+4)+3]$, $pP_{k+1}$ = $\varphi([p_{k}+1,$ $p_{k}(p_k+4)+3])$ = the set of all primes in $IpP_{k+1}$, $\pi(pP_{k+1})$ = $|pP_{k+1}|$  and $(1-\frac{1}{p_k^z})(1-\frac{1}{p_{k-1}^z})\ldots(1-\frac{1}{3^z})(1-\frac{1}{2^z})\zeta(z)$ = $1+\frac{1}{n_{1}^z}+\frac{1}{n_{2}^z}+\ldots$. Then $n_1$ = $p_{k+1}$, $n_2$ = $p_{k+2}$, $\ldots$, $n_{\pi(pP_{k+1})}$ = $p_{k+\pi(pP_{k+1})}\in pP_{k+1}$ = $\varphi([p_{k}+1,$ $p_{k}(p_k+4)+3])$. }  
\end{prop}
\begin{proof} For $Re(z) > 1$, we have

$\zeta(z)$ = $\frac{1}{1^z}+\frac{1}{2^z}+\frac{1}{3^z}+\frac{1}{4^z}+\ldots$.
\\
This implies, 

$\zeta(z)$ = $1+\frac{1}{2^z}(1+\frac{1}{2^z}+\frac{1}{3^z}+\frac{1}{4^z}+\ldots)$ $+\frac{1}{3^z}+\frac{1}{5^z}+\frac{1}{7^z}+\frac{1}{9^z}+\ldots$.
\\
$\Rightarrow$ $(1-\frac{1}{2^z})\zeta(z)$ = $1+\frac{1}{3^z}+\frac{1}{5^z}+\frac{1}{7^z}+\frac{1}{9^z}+\ldots$.

\hspace{1cm} = $1+\frac{1}{2^z}+\frac{1}{3^z}+\frac{1}{4^z}+\ldots$
- $\frac{1}{2^z}(1+\frac{1}{2^z}+\frac{1}{3^z}+\frac{1}{4^z}+\ldots)$.

\hspace{1cm} = 1 + sum of the terms of $\zeta(z)$ without multiples of $\frac{1}{2^z}$.
\\
Similarly, we get,

$(1-\frac{1}{3^z})(1-\frac{1}{2^z})\zeta(z)$ = $1+\frac{1}{5^z}+\frac{1}{7^z}+\frac{1}{11^z}+\frac{1}{13^z}+\frac{1}{17^z}+\frac{1}{19^z}+\frac{1}{23^z}+\ldots$.

\hspace{1cm} = 1 + sum of the terms of $\zeta(z)$ without multiples of $\frac{1}{2^z}$ or $\frac{1}{3^z}$.

Continuing the above process and using Theorem \ref{c7}$(d)$, we get the result. 
\end{proof}

\begin{rem} Number of primes in a pocket of primes plays an important role in the process of generating primes from successive pockets of primes and also measure the power/speed of the algorithm that generates primes. Legendre \cite{l59} obtained a formula to calculate $\pi(x)$, the number of primes $\leq x$, $x\in\mathbb{N}$. Based on this, we obtain, the number of primes in $[p_k+1, p_k(p_k+4)+3]$,

 $\pi([p_k+1, p_k(p_k+4)+3]) = \pi(p_k(p_k+4)+3)-\pi(p_k)$ 

\hspace{1.4cm} $= \pi(p_k(p_k+4)+3)-k$, $k\in\mathbb{N}$.
\end{rem}
Thus, we get the following lemma.

\begin{lemma} \label{l5} {\rm Let $x$ = $p_k$ and $y$ = $p_k(p_k+4)+3$, $k\in\mathbb{N}$. Then, with usual notation, $\pi([p_k+1, p_k(p_k+4)+3])$ = $\pi([x, y])$ = $\pi(p_k(p_k+4)+3)-k$ 

\hfill = $y-x -$ $(\sum_{i=1}^k {(\left\lfloor \frac{y}{p_i} \right\rfloor - \left\lfloor \frac{x}{p_i} \right\rfloor)})$ 
+  $(\sum_{i<j}^k {(\left\lfloor \frac{y}{p_ip_j} \right\rfloor - \left\lfloor \frac{x}{p_ip_j} \right\rfloor)})$ - $\ldots$. }
\end{lemma}
\begin{proof} $\pi([p_k+1, p_k(p_k+4)+3]) =$ $\pi([x, y])$ $= \pi(y)-\pi(x)$ $= \pi(y)-k$.
\\
Also, $\pi([p_k+1, p_k(p_k+4)+3]) =$

\hspace{1cm} $= |[p_k+1, p_k(p_k+4)+3]| - |\mathbb{C}([p_k+1, p_k(p_k+4)+3])|$.

\hspace{1cm} $= |[x+1, y]| - |\mathbb{C}([x+1, y])| = y-x -|\mathbb{C}([x+1, y])|$.

\hspace{1cm} $= (p_k^2+4p_k+3-p_k)$ $-$ $(\sum_{i=1}^k {(\left\lfloor \frac{y}{p_i} \right\rfloor - \left\lfloor \frac{x}{p_i} \right\rfloor)})$ 
\begin{center}
$+$  $(\sum_{i<j}^k {(\left\lfloor \frac{y}{p_ip_j} \right\rfloor - \left\lfloor \frac{x}{p_ip_j} \right\rfloor)})$ $-$ $\ldots$.
\end{center}
\hspace{.2cm} $= (y-x)$ $-$ $(\sum_{i=1}^k {(\left\lfloor \frac{y}{p_i} \right\rfloor - \left\lfloor \frac{x}{p_i} \right\rfloor)})$ $+$  $(\sum_{i<j}^k {(\left\lfloor \frac{y}{p_ip_j} \right\rfloor - \left\lfloor \frac{x}{p_ip_j} \right\rfloor)})$ $-$ $\ldots$.
\end{proof} 
\begin{rem} Lemma \ref{l5} provides a formula for $\pi([p_k+1, p_k(p_k+4)+3])$. But it is very difficult when $k$ is large and through computer programs using the algorithm given for the third method we can calculate $\pi([p_k+1, p_k(p_k+4)+3])$ for different values of $k$.   
\end{rem}
\noindent
\textbf{Conclusion.}
\begin{enumerate} 
\item One can study different sequences of orders of pockets of primes that cover the set of all primes, $\mathbb{P}$.
\item One can study the representation of each natural number by its prime factors instead of successors.
\end{enumerate}   
\noindent
%Table-1
\begin{table}
\caption{$\pi(pP_j)$, order of prime pockets in the three Methods}
\begin{center}
\scalebox{.99}{
\begin{tabular}{||c||*{10}{c|c|c|c|c|c|c|c||}}\hline \hline
$\pi(pP_j)=n_j$ &  $n_1$ & $n_2$ & $n_3$ & $n_4$ & $n_5$ & $n_6$ & $n_7$ & $n_8$  
\\ \hline \hline
Method $1$ & $1$ & $1$ & $1$ & $1$ & $2$ & $3$ & $5$ &  $9$ 
\\ \hline \hline
Method $2$ & $1$ & $1$ & $2$  & $11$ & $314$ & $339730$ &  $- -$ & $- -$ 
\\ \hline \hline
Method $3$ & $2$ & $7$ & $105$  & $32622$ & $- -$ & $- -$ &  $- -$ & $- -$ 
\\
\hline \hline
\end{tabular}}
\end{center}
\end{table} 

\textbf{Acknowledgement.}\quad I express my sincere thanks to the Central University of Kerala, Kasaragod, Kerala, India and St. Jude's College, Thoothoor, Kanyakumari District, Tamil Nadu, India for providing facilities to do this research work. I also express my sincere thanks to Prof. S. Krishnan (late) and Dr. V. Mohan, Department of Mathematics, Thiyagarayar College of Engg., Madurai, Tamil Nadu, India and Prof. M.I. Jinnah (late), Department of Mathematics, University of Kerala, Thiruvanandapuram, Kerala, India for their support and encouragements to do research.

\end{document}